\documentclass[runningheads]{llncs}
\usepackage{graphicx}
\usepackage{amsmath,amssymb}
\usepackage{comment}
\usepackage{algpseudocode}
\usepackage{algorithm}
\usepackage{tikz}
\usepackage{pgfplots}
\usetikzlibrary{shapes,snakes}

\errorcontextlines\maxdimen

\makeatletter
\newcommand*{\algrule}[1][\algorithmicindent]{\makebox[#1][l]{\hspace*{.5em}\thealgruleextra\vrule height \thealgruleheight depth \thealgruledepth}}%
\newcommand*{\thealgruleextra}{}
\newcommand*{\thealgruleheight}{.75\baselineskip}
\newcommand*{\thealgruledepth}{.25\baselineskip}

\newcount\ALG@printindent@tempcnta
\def\ALG@printindent{%
	\ifnum \theALG@nested>0
	\ifx\ALG@text\ALG@x@notext
	\else
	\unskip
	\addvspace{-1pt}
	\ALG@printindent@tempcnta=1
	\loop
	\algrule[\csname ALG@ind@\the\ALG@printindent@tempcnta\endcsname]%
	\advance \ALG@printindent@tempcnta 1
	\ifnum \ALG@printindent@tempcnta<\numexpr\theALG@nested+1\relax
	\repeat
	\fi
	\fi
}%
\usepackage{etoolbox}
\patchcmd{\ALG@doentity}{\noindent\hskip\ALG@tlm}{\ALG@printindent}{}{\errmessage{failed to patch}}
\makeatother

\newbox\statebox
\newcommand{\myState}[1]{%
	\setbox\statebox=\vbox{#1}%
	\edef\thealgruleheight{\dimexpr \the\ht\statebox+1pt\relax}%
	\edef\thealgruledepth{\dimexpr \the\dp\statebox+1pt\relax}%
	\ifdim\thealgruleheight<.75\baselineskip
	\def\thealgruleheight{\dimexpr .75\baselineskip+1pt\relax}%
	\fi
	\ifdim\thealgruledepth<.25\baselineskip
	\def\thealgruledepth{\dimexpr .25\baselineskip+1pt\relax}%
	\fi
	\State #1%
	\def\thealgruleheight{\dimexpr .75\baselineskip+1pt\relax}%
	\def\thealgruledepth{\dimexpr .25\baselineskip+1pt\relax}%
}

\begin{document}
\title{An iterative ILP approach for constructing a Hamiltonian decomposition of a regular multigraph \thanks{This research was supported by P.G. Demidov Yaroslavl State University Project VIP-016}}
\titlerunning{An iterative ILP approach for constructing a Hamiltonian decomposition}
%
\author{Andrey Kostenko\orcidID{0000-0001-7482-9821} \and
Andrei Nikolaev\orcidID{0000-0003-4705-2409}}
\authorrunning{A. Kostenko, A. Nikolaev}
%
\institute{P.\,G. Demidov Yaroslavl State University, Yaroslavl, Russia,
\email{andreykostenko99@gmail.com, andrei.v.nikolaev@gmail.com}}
\maketitle              
\begin{abstract}
A Hamiltonian decomposition of a regular graph is a partition of its edge set into Hamiltonian cycles. The problem of finding edge-disjoint Hamiltonian cycles in a given regular graph has many applications in combinatorial optimization and operations research. Our motivation for this problem comes from the field of polyhedral combinatorics, as a sufficient condition for vertex nonadjacency in the 1-skeleton of the traveling salesperson polytope can be formulated as the Hamiltonian decomposition problem in a 4-regular multigraph with one forbidden decomposition. 

In our approach, the algorithm starts by solving the relaxed 2-matching problem, then iteratively generates subtour elimination constraints for all subtours in the solution and solves the corresponding ILP-model to optimality. The procedure is enhanced by the local search heuristic based on chain edge fixing and cycle merging operations. In the computational experiments, the iterative ILP algorithm showed comparable results with the previously known heuristics on undirected multigraphs and significantly better performance on directed multigraphs.


\end{abstract}

\section{Introduction}

\textit{A Hamiltonian decomposition} of a regular graph is a partition of its edge set into Hamiltonian cycles. The problem of finding edge-disjoint Hamiltonian cycles in a given regular graph has different applications in combinatorial optimization \cite{krar:1995}, coding theory \cite{bae:bose:2003,bail:2009}, privacy-preserving distributed mining algorithms \cite{clif:kant:vaid:2002}, analysis of interconnection networks \cite{hung:2011} and other areas. 
Our motivation for this problem comes from the field of polyhedral combinatorics.

We consider a classic traveling salesperson problem: given a complete weighted graph (or digraph) $K_n=(V,E)$, it is required to find a Hamiltonian cycle of minimum weight.
We denote by $HC_{n}$ the set of all Hamiltonian cycles in $K_{n}$.
With each Hamiltonian cycle $x \in HC_n$ we associate a characteristic vector $x^v \in \mathbb{R}^{E}$ by the following rule:
\[
x^v_e = 
\begin{cases}
1,& \text{ if the cycle } x \text{ contains an edge } e \in E,\\
0,& \text{ otherwise. }
\end{cases}
\]
The polytope
\[\operatorname{TSP}(n) = \operatorname{conv} \{x^v \ | \ x \in HC_n \}\]
is called \textit{the symmetric traveling salesperson polytope}.

\textit{The asymmetric traveling salesperson polytope} $\operatorname{ATSP}(n)$ is defined similarly as the convex hull of characteristic vectors of all possible Hamiltonian cycles in the complete digraph $K_{n}$.

The 1-\textit{skeleton} of a polytope $P$ is the graph whose vertex set is the vertex set of $P$ and edge set is the set of one-dimensional faces of $P$. 
The study of 1-skeleton is of interest, since, on the one hand, some algorithms for perfect matching, set covering, independent set, a ranking of objects, and problems with fuzzy measures are based on the vertex adjacency relation in 1-skeleton and the local search technique (see, for example, \cite{ag:katz:tol:2017,bali:1985,cher:ham:1987,comb:mir:2010}). On the other hand, such characteristics of 1-skeleton as the diameter and clique number, estimate the time complexity for different computation models and classes of algorithms \cite{bond:1983,bond:nik:2013,grot:padb:1985}.

Unfortunately, the classic result by Papadimitriou stands in the way of studying the 1-skeleton of the traveling salesperson polytope.

\begin{theorem} [Papadimitriou \cite{papa:1978}]\label{theorem:papa}
	The question of whether two vertices of the polytopes $\operatorname{TSP}(n)$ or $\operatorname{ATSP}(n)$ are nonadjacent is NP-complete. 
\end{theorem}

\section{Formulation of the problem}

We consider a sufficient condition for vertex nonadjacency in 1-skeleton of the traveling salesperson polytope by Rao \cite{rao:1976}.
Let $x = (V,E(x))$ and $y=(V,E(y))$ be two Hamiltonian cycles on the vertex set $V$.
We denote by $x \cup y$ a multigraph $(V,E(x) \cup E(y))$ that contains all edges of both cycles $x$ and $y$.

\begin{lemma}[Rao \cite{rao:1976}]\label{lemma_sufficient}
	Given two Hamiltonian cycles $x$ and $y$, if the multigraph $x \cup y$ contains a Hamiltonian decomposition into edge-disjoint cycles $z$ and $w$ different from $x$ and $y$, then the corresponding vertices $x^v$ and $y^v$ of the polytope $\operatorname{TSP}(n)$ (or $\operatorname{ATSP}(n)$) are not adjacent.
\end{lemma}

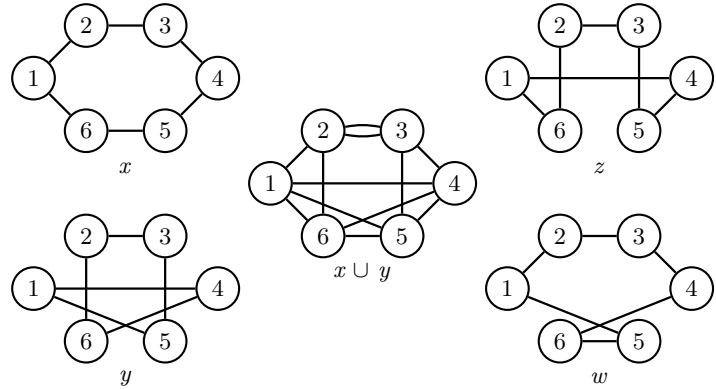
\begin{figure}[t]
	\centering
	\begin{tikzpicture}[scale=0.7]
	\begin{scope}[every node/.style={circle,thick,draw}]
	\node (1) at (0,0) {1};
	\node (2) at (1,1) {2};
	\node (3) at (2.5,1) {3};
	\node (4) at (3.5,0) {4};
	\node (5) at (2.5,-1) {5};
	\node (6) at (1,-1) {6};
	\end{scope}
	\draw [line width=0.3mm] (1) edge (2);
	\draw [line width=0.3mm] (2) edge (3);
	\draw [line width=0.3mm] (3) edge (4);
	\draw [line width=0.3mm] (4) edge (5);
	\draw [line width=0.3mm] (5) edge (6);
	\draw [line width=0.3mm] (6) edge (1);
	\draw (1.75, -1.7) node{\textit{x}};
	
	\begin{scope}[yshift=-4cm]
	\begin{scope}[every node/.style={circle,thick,draw}]
	\node (1) at (0,0) {1};
	\node (2) at (1,1) {2};
	\node (3) at (2.5,1) {3};
	\node (4) at (3.5,0) {4};
	\node (5) at (2.5,-1) {5};
	\node (6) at (1,-1) {6};
	\end{scope}
	\draw [line width=0.3mm] (1) edge (4);
	\draw [line width=0.3mm] (4) edge (6);
	\draw [line width=0.3mm] (6) edge (2);
	\draw [line width=0.3mm] (2) edge (3);
	\draw [line width=0.3mm] (3) edge (5);
	\draw [line width=0.3mm] (5) edge (1);
	\draw (1.75, -1.7) node{\textit{y}};
	\end{scope}
	
	\begin{scope}[xshift=4.5cm,yshift=-2cm]
	\begin{scope}[every node/.style={circle,thick,draw}]
	\node (1) at (0,0) {1};
	\node (2) at (1,1) {2};
	\node (3) at (2.5,1) {3};
	\node (4) at (3.5,0) {4};
	\node (5) at (2.5,-1) {5};
	\node (6) at (1,-1) {6};
	\end{scope}
	\draw [line width=0.3mm] (1) edge (2);
	\draw [line width=0.3mm, bend right=10] (2) edge (3);
	\draw [line width=0.3mm] (3) edge (4);
	\draw [line width=0.3mm] (4) edge (5);
	\draw [line width=0.3mm] (5) edge (6);
	\draw [line width=0.3mm] (6) edge (1);
	\draw [line width=0.3mm] (1) edge (4);
	\draw [line width=0.3mm] (4) edge (6);
	\draw [line width=0.3mm] (6) edge (2);
	\draw [line width=0.3mm, bend left=10] (2) edge (3);
	\draw [line width=0.3mm] (3) edge (5);
	\draw [line width=0.3mm] (5) edge (1);	
	\draw (1.75, -1.7) node{\textit{x $\cup$ y}};
	\end{scope}

	\begin{scope}[xshift=9cm]
	\begin{scope}[every node/.style={circle,thick,draw}]
	\node (1) at (0,0) {1};
	\node (2) at (1,1) {2};
	\node (3) at (2.5,1) {3};
	\node (4) at (3.5,0) {4};
	\node (5) at (2.5,-1) {5};
	\node (6) at (1,-1) {6};
	\end{scope}
	\draw [line width=0.3mm] (1) edge (4);
	\draw [line width=0.3mm] (4) edge (5);
	\draw [line width=0.3mm] (5) edge (3);
	\draw [line width=0.3mm] (3) edge (2);
	\draw [line width=0.3mm] (2) edge (6);
	\draw [line width=0.3mm] (6) edge (1);
	\draw (1.75, -1.7) node{\textit{z}};
	
	\begin{scope}[yshift=-4cm]
	\begin{scope}[every node/.style={circle,thick,draw}]
	\node (1) at (0,0) {1};
	\node (2) at (1,1) {2};
	\node (3) at (2.5,1) {3};
	\node (4) at (3.5,0) {4};
	\node (5) at (2.5,-1) {5};
	\node (6) at (1,-1) {6};
	\end{scope}
	\draw [line width=0.3mm] (1) edge (2);
	\draw [line width=0.3mm] (2) edge (3);
	\draw [line width=0.3mm] (3) edge (4);
	\draw [line width=0.3mm] (4) edge (6);
	\draw [line width=0.3mm] (6) edge (5);
	\draw [line width=0.3mm] (5) edge (1);
	\draw (1.75, -1.7) node{\textit{w}};
	\end{scope}
	\end{scope}
	\end{tikzpicture}
	\caption{The multigraph $x \cup y$ has two different Hamiltonian decompositions}
	\label{image:not_adjacent}
\end{figure}

From a geometric point of view, the sufficient condition means that the segment connecting two vertices $x^v$ and $y^v$ intersects with the segment connecting two other vertices $z^v$ and $w^v$ of the polytope $\operatorname{TSP}(n)$ (or $\operatorname{ATSP}(n)$ correspondingly). Thus, the vertices $x^v$ and $y^v$ cannot be adjacent. An example of a satisfied sufficient condition is shown in Fig.~\ref{image:not_adjacent} (see also~\cite{nik:kozl:2020}).

We formulate the sufficient condition for vertex nonadjacency of the traveling salesperson polytope as a combinatorial problem.

\vspace{1mm}

\textbf{Hamiltonian decomposition with one forbidden decomposition.}

\textsc{Instance.} Let $x$ and $y$ be two Hamiltonian cycles.

\textsc{Question.} Does the multigraph $x \cup y$ contain a pair of edge-disjoint Hamiltonian cycles $z$ and $w$ different from $x$ and $y$?

\vspace{1mm}

Thus, we consider a version of a Hamiltonian decomposition problem of a special form.
By construction, the union multigraph $x \cup y$ always contains the Hamiltonian decomposition into $x$ and $y$, but in our case, this decomposition is forbidden, and we need to find another decomposition if it exists. 

Note that testing of whether a graph has a Hamiltonian decomposition is NP-complete, even for 4-regular undirected graphs and 2-regular directed graphs \cite{per:1984}.
Therefore, instead of Rao's sufficient condition, various polynomially solvable special cases of the vertex nonadjacency problem have been studied in the literature.
In particular, the polynomial sufficient conditions for the pyramidal tours \cite{bond:nik:2018}, pyramidal tours with step-backs \cite{nik:2019}, and pedigrees \cite{arth:2006,arth:2013} are known.
However, all of them are weaker than the sufficient condition by Rao. 

The Hamiltonian decomposition problem of the considered form was introduced in \cite{kozl:nik:2019} and later studied in \cite{nik:kozl:2020}, where two heuristic algorithms were proposed.
The set of feasible solutions in both algorithms consists of all possible decompositions of the multigraph $x \cup y$ into edge-disjoint 2-factors $z$ and $w$.
Recall that \textit{a $2$-factor} (or \textit{a perfect 2-matching}) of a graph $G$ is a subset of edges of $G$ such that every vertex is incident with exactly two edges.
The differences are as follows:
\begin{itemize}
	\item the simulated annealing algorithm from \cite{kozl:nik:2019} repeatedly finds 2-factors $z$ and $w$ through the reduction to random perfect matching \cite{tutt:1954} 
	until it obtains a pair of Hamiltonian cycles different from $x$ and $w$;
	\item the general variable neighborhood search algorithm from \cite{nik:kozl:2020} adds to the previous algorithm several neighborhood structures and cycle merging operations combined in the basic variable neighborhood
	descent approach \cite{duar:san:mlad:tod:2018}.
\end{itemize}
Heuristic algorithms have proven to be very efficient on instances with an existing solution, especially on undirected graphs. However, on instances without a solution, the heuristics face significant difficulties.

In this paper, we propose two exact algorithms for solving the Hamiltonian decomposition problem with a forbidden decomposition and verifying vertex nonadjacency of the traveling salesperson polytope.
The first algorithm iteratively generates ILP models for the problem.
The second algorithm combines the first one and the modified local search heuristic from \cite{nik:kozl:2020}.

\section{Iterative integer linear programming} \label{sec:iterativeILP}

Let $x = (V,E(x))$, $y = (V,E(y))$, $x \cup y = (V,E = E(x) \cup E(y))$.
With each edge $e \in E$ we associate the variable
\[
x_e = \begin{cases}
1,& \text{if } e \in z,\\
0,& \text{if } e \in w.
\end{cases}
\]

We adapt the classic ILP formulation of the traveling salesperson problem by Dantzig, Fulkerson and Johnson \cite{dantz:falk:john:1954} into the following ILP model for the considered Hamiltonian decomposition problem:
\begin{align}
&\sum_{e \in E} x_e = |V|,  \label{ILP_tour}\\
&\sum_{e \in E_v} x_e = 2, 		&\forall v \in V, \label{ILP_vertex_degree}\\
&\sum_{e \in E(x) \backslash E(y)} x_e \leq |V|-|E(x) \cap E(y)|-2,  	\label{ILP_not_x}\\
&\sum_{e \in E(y) \backslash E(x)} x_e \leq |V|-|E(x) \cap E(y)|-2,	\label{ILP_not_y}\\
&\sum_{e \in E_S} x_e \leq |S|-1, 		&\forall S \subset V,	\label{ILP_subtour_elim_1}\\
&\sum_{e \in E_S} x_e \geq |E_S| - |S| + 1,		&\forall S \subset V,	\label{ILP_subtour_elim_0}\\
&x_e \in \{0,1\}, 		&\forall e \in E. 	\label{ILP_variables}
\end{align}

In the following, we elaborate on the model.
The multigraph $x \cup y$ contains $2 |V|$ edges. The constraint (\ref{ILP_tour}) guarantees that both components $z$ ($x_e = 1$) and $w$ ($x_e = 0$) receive exactly $|V|$ edges.

We denote by $E_v$ the set of all edges incident to the vertex $v$ in $x \cup y$. By the \textit{vertex degree constraint} (\ref {ILP_vertex_degree}) the degree of each vertex in $z$ and $w$ is equal to 2.

The vertex degree constraint for directed graphs is slightly different. Let for some vertex $v \in V $: $e_1$ and $e_2$ be two incoming edges, and $u_1$ and $u_2$ be two outgoing edges. Then the constraints (\ref {ILP_vertex_degree}) will take the form:
\begin{align*}
x_{e_1} + x_{e_2} = 1,\\ 
x_{u_1} + x_{u_2} = 1,
\end{align*}
with exactly one incoming edge and outgoing edge for each vertex in the solution.

The constraints (\ref{ILP_not_x})-(\ref{ILP_not_y}) forbid the Hamiltonian cycles $x$ and $y$ as a solution.
If we consider a general Hamiltonian decomposition problem of a 4-regular multigraph without reference to the vertex adjacency in the traveling salesperson polytope, then these constraints can be omitted.

Finally, the inequalities (\ref{ILP_subtour_elim_1})-(\ref{ILP_subtour_elim_0}) are known as the \textit{subtour elimination constraints} (\textit{SEC}), which forbid solutions consisting of several disconnected tours.
Here $S$ is a subset of $V$, $E_S$ is the set of all edges from $E$ with both vertices belonging to $S$:
\[E_S = \{(u, v) \in E: \ u, v \in S \}. \]

The main problem with the subtour elimination constraints is that there are exponentially many of them: two for each subset of $S \subset V$, i.e. $\Omega (2^{|V|})$.
Therefore, the idea of the first algorithm is as follows. 
We start with the relaxed model (\ref{ILP_tour})-(\ref{ILP_not_y}),(\ref{ILP_variables}) of the basic 2-matching problem with $O(V)$ constraints. By ILP-solver we obtain an integer point that corresponds to the pair of 2-factors $z$ and $w$.
Then we find all subtours in $z$ and $w$, add the corresponding subtour elimination constraints (\ref{ILP_subtour_elim_1}) and (\ref{ILP_subtour_elim_0}) into the model, and iteratively repeat this procedure.
The algorithm stops either by finding the Hamiltonian decomposition into cycles $z$ and $w$ or by obtaining an infeasible model that does not contain any integer points.
This approach is inspired by the algorithm for the traveling salesperson problem from \cite{pfer:stan:2017} and is summarized in Algorithm~\ref{Alg:ILP}.

\begin{algorithm}[t]
	\caption{Iterative integer linear programming algorithm}
	\label{Alg:ILP}
	\begin{algorithmic}[0]
		\Procedure{IterativeILP}{$x \cup y$}
		\State Define current model as (\ref{ILP_tour})-(\ref{ILP_not_y}),(\ref{ILP_variables})		\Comment{relaxed 2-matching problem}
		\While {the model is feasible} 
			\State $z,w \gets$ an integer point of the current model by an ILP-solver
			\If{$z$ and $w$ is a Hamiltonian decomposition}
				\State \Return Hamiltonian decomposition $z$ and $w$
			\EndIf
			\State For all subtours in $z$ and $w$ add the SEC (\ref{ILP_subtour_elim_1}) and (\ref{ILP_subtour_elim_0}) into the model	
		\EndWhile	
		\State \Return Hamiltonian decomposition does not exist
		\EndProcedure		
	\end{algorithmic}
\end{algorithm}

\section{Local search} \label{sec:local_search}

To improve the performance, we enhance Algorithm~\ref{Alg:ILP} with the local search heuristic.
The neighborhood structure is a significantly modified version of the first neighborhood in the GVNS algorithm \cite{nik:kozl:2020} with a new recursive chain edge fixing procedure.

\subsection{Feasible set}

\begin{figure}[p]
	\centering
	\begin{tikzpicture}[scale=0.85]
	\begin{scope}[every node/.style={circle,thick,draw}]
	\node (a1) at (0,0) {1};
	\node (a2) at (1,1) {2};
	\node (a3) at (2.5,1) {3};
	\node (a4) at (3.5,0) {4};
	\node (a5) at (2.5,-1) {5};
	\node (a6) at (1,-1) {6};
	\end{scope}
	
	\node at (1.75, -1.75) {$x \cup y$};
	
	\draw [thick] (a1) edge (a2);
	\draw [thick] (a1) edge (a6);
	\draw [thick] (a2) edge (a6);
	\draw [thick] (a3) edge (a4);
	\draw [thick] (a3) edge (a5);
	\draw [thick] (a4) edge (a5);
	
	\draw [thick, bend left=25] (a1) edge (a2);
	\draw [thick] (a1) edge (a3);
	\draw [thick] (a2) edge (a3);
	\draw [thick] (a4) edge (a6);
	\draw [thick] (a5) edge (a6);
	\draw [thick, bend left=25] (a4) edge (a5);

	\begin{scope}[xshift=-4.75cm]
	\begin{scope}[every node/.style={circle,thick,draw}]
	\node (b1) at (0,0) {1};
	\node (b2) at (1,1) {2};
	\node (b3) at (2.5,1) {3};
	\node (b4) at (3.5,0) {4};
	\node (b5) at (2.5,-1) {5};
	\node (b6) at (1,-1) {6};
	\end{scope}
	
	\draw [thick] (b1) edge (b2);
	\draw [thick] (b1) edge (b6);
	\draw [thick] (b2) edge (b6);
	\draw [thick] (b3) edge (b4);
	\draw [thick] (b3) edge (b5);
	\draw [thick] (b4) edge (b5);
	\node at (1.75, -1.75) {$z$};

	\end{scope}

	\begin{scope}[xshift=4.75cm]
	\begin{scope}[every node/.style={circle,thick,draw}]
	\node (c1) at (0,0) {1};
	\node (c2) at (1,1) {2};
	\node (c3) at (2.5,1) {3};
	\node (c4) at (3.5,0) {4};
	\node (c5) at (2.5,-1) {5};
	\node (c6) at (1,-1) {6};
	\end{scope}
	
	\node at (1.75, -1.75) {$w = (x \cup y) \backslash z$};
	
	\draw [thick, bend left=25] (c1) edge (c2);
	\draw [thick] (c1) edge (c3);
	\draw [thick] (c2) edge (c3);
	\draw [thick] (c4) edge (c6);
	\draw [thick] (c5) edge (c6);
	\draw [thick, bend left=25] (c4) edge (c5);	
	\end{scope}
	\end{tikzpicture}
	\caption{The multigraph $x \cup y$ and its two edge-disjoint 2-factors}
	\label{image:2-factors}
\end{figure}
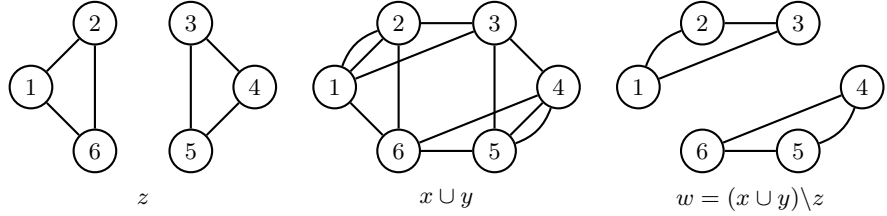

Every solution of the ILP model (\ref{ILP_tour})-(\ref{ILP_not_y}),(\ref{ILP_variables}) with partial subtour elimination constraints corresponds to the pair $z$ and $w$ of edge-disjoint 2-factors of the multigraph $x \cup y$ (Fig.~\ref{image:2-factors}).
We compose a set of feasible solutions for the local search algorithm from all possible pairs of edge-disjoint 2-factors of the multigraph $x \cup y$.

\subsection{Objective function}

As the objective function to minimize, we choose the total number of connected components in the 2-factors $z$ and $w$.
If it equals 2, then $z$ and $w$ are Hamiltonian cycles.

\subsection{Neighborhood structure for directed graphs}

The main difference between the neighborhood structure presented in this chapter and those described in \cite{nik:kozl:2020} is the \textit{chain edge fixing procedure}.
We divide the edges of $z$ and $w$ into two classes:
\begin{itemize}
	\item \textit{unfixed} edges that can be moved between $z$ and $w$ to get a neighboring solution;
	\item edges that are \textit{fixed} in $z$ or $w$ and cannot be moved.
\end{itemize}

\begin{figure}[p]
	\centering
	\begin{tikzpicture}[scale=0.9]
	\begin{scope}[every node/.style={circle,thick,draw}]
	\node (i) at (0,0) {$i$};
	\node (j) at (3,0) {$j$};
	\node (k) at (3,1.5) {$k$};
	\node (h) at (0,-1.5) {$h$};
	\end{scope}
	
	\draw [->,>=stealth,thick] (i) edge node[below]{$z$} (j);
	\draw [->,>=stealth,thick,dashed] (i) edge node[above]{$w$} (k);
	\draw [->,>=stealth,thick,dashed] (h) edge node[below]{$w$} (j);
	
	\end{tikzpicture}
	\caption{Fixing the edge $(i,j)$ in $z$}
	\label{Fig_fixed_edges}
\end{figure}
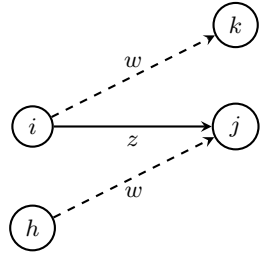

The idea is that one fixed edge starts a recursive chain of fixing other edges.
For example, we consider a directed 2-regular multigraph $x \cup y$ with all indegrees and outdegrees are equal to 2. Let us choose some edge $(i,j)$ and fix it in the component $z$, then the second edge $(i,k)$ outgoing from $i$ and the second edge $(h,j)$ incoming into $j$ obviously cannot get into $z$. We will fix these edges in $w$ (Fig.~\ref{Fig_fixed_edges}). 
In turn, the edges $(i,k)$ and $(h,j)$, fixed in $w$, start the recursive chains of fixing edges in $z$, etc.

At the initial step, we fix in $z$ and $w$ a copy of each multiple edge of the $x \cup y$, since both copies obviously cannot end up in the same Hamiltonian cycle.
We construct a neighboring solution by choosing an edge of $z$, moving it to $w$, and running the chain edge fixing procedure to restore the correct 2-factors.
If the number of connected components in $z$ and $w$ has decreased, we find all subtours in $z$ and $w$, add the corresponding subtour elimination constraints (\ref{ILP_subtour_elim_1}) and (\ref{ILP_subtour_elim_0}) into the model, proceed to a new solution and restart the local search.
This procedure is summarized in Algorithm~\ref{Alg:LS_directed}.

\begin{algorithm}[p]
	\caption{Local search for directed graphs}\label{Alg:LS_directed}
	\begin{algorithmic}[0]
		\Procedure{Chain\_Edge\_Fixing\_Directed}{$(i,j)$ in $z$}		
		\State Fix the edge $(i,j)$ in $z$ and mark it as checked to avoid double checking
		\If{the edge $(i,k)$ is not fixed}
			\State \Call{Chain\_Edge\_Fixing\_Directed}{$(i,k)$ in $w$}
		\EndIf
		\If{the edge $(h,j)$ is not fixed}
			\State \Call{Chain\_Edge\_Fixing\_Directed}{$(h,j)$ in $w$}
		\EndIf
		\EndProcedure
		\vspace{2mm}
		\Procedure{Local\_Search\_Directed}{$z,w$}
			\State Fix the multiple edges in $z$ and $w$
			\Repeat 
				\State Shuffle the edges of $z$ in random order	
				\For {each unchecked and unfixed edge $(i,j)$ in $z$}
					\State \Call {Chain\_Edge\_Fixing\_Directed} {$(i,j)$ in $w$} 	\Comment{Move $(i,j)$ from $z$ to $w$}
					\If {the number of connected components in $z$ and $w$ has decreased}
						\State For all subtours in $z$ and $w$ add the SEC (\ref{ILP_subtour_elim_1}) and (\ref{ILP_subtour_elim_0}) into the model
						\State Proceed to a new solution and restart the local search
					\EndIf
					\State Restore the original $z$ and $w$ and unfix all non-multiple edges
				\EndFor
			\Until  all edges of $z$ are checked and no improvement found			\Comment{A local minimum}
			\State \Return $z$ and $w$
		\EndProcedure
	\end{algorithmic}
\end{algorithm}

Note that although the chain edge fixing procedure at each step can call up to two of its recursive copies, the total complexity is linear ($O(V)$), since each edge can only be fixed once, and $|E| = 2 |V|$.
Thus, the size of the neighborhood is equal to the number of edges in $z$, i.e. $O (V)$, and the total complexity of exploring the neighborhood is quadratic $O(V^2)$.

\subsection{Neighborhood structure for undirected graphs}

The neighborhood structure for undirected graphs is similar, we choose an edge of $z$, move it to $w$ and run the chain edge fixing procedure.

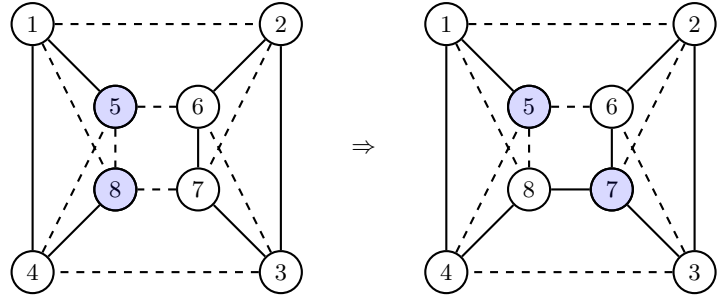
\begin{figure}[p]
\centering
\begin{tikzpicture}[scale=1.1]	
	\begin{scope}[every node/.style={circle,draw,thick,inner sep=3pt}]
	\node (1) at (0,3) {$1$};
	\node (2) at (3,3) {$2$};
	\node (3) at (3,0) {$3$};
	\node (4) at (0,0) {$4$};
	\node (5) at (1,2) {$5$};
	\node (6) at (2,2) {$6$};
	\node (7) at (2,1) {$7$};
	\node (8) at (1,1) {$8$};
	\end{scope}
	
	\node [circle,draw,thick,inner sep=3pt,fill=blue,fill opacity=0.15] (5) at (1,2) {$5$};
	\node [circle,draw,thick,inner sep=3pt,fill=blue,fill opacity=0.15] (8) at (1,1) {$8$};
	
	\draw [thick] (1) -- (5);
	\draw [thick] (8) -- (4);
	\draw [thick] (4) -- (1);
	
	\draw [thick] (2) -- (3);
	\draw [thick] (3) -- (7);
	\draw [thick] (6) -- (7);
	\draw [thick] (2) -- (6);

	\draw [thick,dashed] (1) -- (2);
	\draw [thick,dashed] (2) -- (7);
	\draw [thick,dashed] (7) -- (8);
	\draw [thick,dashed] (1) -- (8);
	
	\draw [thick,dashed] (3) -- (6);
	\draw [thick,dashed] (5) -- (6);
	\draw [thick,dashed] (4) -- (5);
	\draw [thick,dashed] (3) -- (4);
	
	\draw [thick,dashed] (5) -- (8);
	
	\node at (4,1.5) {$\Rightarrow$};
	
	\begin{scope}[xshift=5cm]
	\begin{scope}[every node/.style={circle,draw,thick,inner sep=3pt}]
	\node (1) at (0,3) {$1$};
	\node (2) at (3,3) {$2$};
	\node (3) at (3,0) {$3$};
	\node (4) at (0,0) {$4$};
	\node (5) at (1,2) {$5$};
	\node (6) at (2,2) {$6$};
	\node (7) at (2,1) {$7$};
	\node (8) at (1,1) {$8$};
	\end{scope}
	
	\node [circle,draw,thick,inner sep=3pt,fill=blue,fill opacity=0.15] (5) at (1,2) {$5$};
	\node [circle,draw,thick,inner sep=3pt,fill=blue,fill opacity=0.15] (7) at (2,1) {$7$};
	
	\draw [thick] (1) -- (5);
	\draw [thick] (8) -- (4);
	\draw [thick] (4) -- (1);
	
	\draw [thick] (2) -- (3);
	\draw [thick] (3) -- (7);
	\draw [thick] (6) -- (7);
	\draw [thick] (2) -- (6);
	
	\draw [thick] (7) -- (8);

	\draw [thick,dashed] (1) -- (2);
	\draw [thick,dashed] (2) -- (7);
	\draw [thick,dashed] (1) -- (8);
	
	\draw [thick,dashed] (3) -- (6);
	\draw [thick,dashed] (5) -- (6);
	\draw [thick,dashed] (4) -- (5);
	\draw [thick,dashed] (3) -- (4);
	
	\draw [thick,dashed] (5) -- (8);
	\end{scope}
-\end{tikzpicture}
\caption{Restoring the broken vertex $8$ (blue colored) by moving the random incident edge $(8,7)$ from $w$ (dashed edges) to $z$ (solid edges).}
\label{Fig_restore_broken_vertex}
\end{figure}

The key difference is that after the exchange of edges and the chain edge fixing procedure, there will remain some \textit{broken vertices} in $z$ and $w$ with a degree not equal to 2. We restore the degree of each broken vertex by moving random unfixed incident edges between components $z$ and $w$ (Fig.~\ref{Fig_restore_broken_vertex}).
This procedure is summarized in Algorithm~\ref{Alg:LS_undirected}.

\begin{algorithm}[p]
	\caption{Local search for undirected graphs}\label{Alg:LS_undirected}
	\begin{algorithmic}[0]
		\Procedure{Chain\_Edge\_Fixing\_Undirected}{$(i,j)$ in $z$}		
			\State Fix the edge $(i,j)$ in $z$
			\If{vertex $i$ in $z$ has two incident fixed edges}
				\State \Call{Chain\_Edge\_Fixing\_Undirected}{$(i,k)$ and $(i,h)$ in $w$}
			\EndIf
			\If{vertex $j$ in $z$ has two incident fixed edges}
				\State \Call{Chain\_Edge\_Fixing\_Undirected}{$(j,k)$ and $(j,h)$ in $w$}
			\EndIf
		\EndProcedure
		\vspace{1mm}
		\Procedure{Local\_Search\_Undirected}{$z,w,attemptLimit$}
			\State Fix the multiple edges in $z$ and $w$
			\Repeat
			\State Shuffle the edges of $z$ in random order		
			\For {each unchecked and unfixed edge $(i,j)$ in $z$}
				\State \Call {Chain\_Edge\_Fixing\_Undirected}{$(i,j)$ in $w$}		\Comment{Move $(i,j)$ from $z$ to $w$}
				\For {$i \gets 1$ \textbf{to} $attemptLimit$}				
					\While {$z$ contains a broken vertex $i$ with degree not equal to 2}
						\If {vertex degree of $i$ is equal to 1} 				
							\State Pick a random unfixed edge $(i,k)$ of $w$	\Comment{One missing edge}
							\State \Call {Chain\_Edge\_Fixing\_Undirected}{$(i,k)$ in $z$} 
						\EndIf 
						\If {vertex degree of $i$ is equal to 3} 				
							\State Pick a random unfixed edge $(i,k)$ of $z$;	\Comment{One extra edge}
							\State \Call {Chain\_Edge\_Fixing\_Undirected}{$(i,k)$ in $w$} 
						\EndIf
					\EndWhile
					\If {the number of connected components has decreased}
						\State For all subtours in $z$ and $w$ add the SEC (\ref{ILP_subtour_elim_1}) and (\ref{ILP_subtour_elim_0}) into the model
						\State Proceed to a new solution and restart the local search
					\EndIf
					\State Restore $z$ and $w$ and unfix all non-multiple edges except $(i,j)$
				\EndFor
				\State Unfix the edge $(i,j)$ and mark it as checked
			\EndFor
			\Until all edges of $z$ are checked and no improvement found			\Comment{A local minimum}
			\State \Return $z$ and $w$
		\EndProcedure
	\end{algorithmic}
\end{algorithm}

Since at each step we pick a random edge to restore a broken vertex, the local search for undirected graphs is a randomized algorithm. Therefore, we run several attempts (parameter $attemptLimit$) while constructing each neighboring solution, i.e. exploring several random branches in the search tree. Thus, the size of the neighborhood is equal to $O (V \cdot attemptLimit)$, and the total complexity of exploring the neighborhood is $O(V^2 \cdot attemptLimit)$.

\subsection{Iterative ILP algorithm with local search}

We add the local search heuristic into the iterative ILP algorithm between iterations to improve the performance on instances with an existing Hamiltonian decomposition. 
If the ILP-solver returns a pair of edge-disjoint 2-factors $z$ and $w$ that are not a Hamiltonian decomposition, then we call the local search to minimize the number of connected components.

Note that every time the local search improves the solution, we modify the model by adding the corresponding subtour elimination constraints for all subtours in $z$ and $w$.
Thus, we implement the memory structure and prohibit the algorithm from returning to feasible solutions that have already been explored.

If the heuristic also fails, we restart the ILP-solver on the modified model, and repeat these steps until a Hamiltonian decomposition is found, or the resulting model is infeasible.
This procedure is summarized in Algorithm~\ref{Alg:ILP_LS}.

\begin{algorithm}[h]
	\caption{Iterative ILP algorithm with local search}
	\label{Alg:ILP_LS}
	\begin{algorithmic}[0]
		\Procedure{IterativeILP+LS}{$x \cup y,attemptLimit$}
			\State define the current model as (\ref{ILP_tour})-(\ref{ILP_not_y}),(\ref{ILP_variables})	\Comment{relaxed 2-matching problem}
			\While {the model is feasible}
				\State $z,w \gets$ an integer point of the current model by ILP-solver
				\If{$z$ and $w$ is a Hamiltonian decomposition}
					\State \Return Hamiltonian decomposition $z$ and $w$
				\EndIf
				\State For all subtours in $z$ and $w$ add the SEC (\ref{ILP_subtour_elim_1}) and (\ref{ILP_subtour_elim_0}) into the model
				\If {the graph is directed}
					\State $z,w \gets$  \Call {Local\_Search\_Directed} {$z,w$};
				\Else
					\State $z,w \gets$  \Call {Local\_Search\_Undirected} {$z,w,attemptLimit$};
				\EndIf
				\If{$z$ and $w$ is a Hamiltonian decomposition}
					\State \Return Hamiltonian decomposition $z$ and $w$
				\EndIf
			\EndWhile
			\State \Return Hamiltonian decomposition does not exist
		\EndProcedure		
	\end{algorithmic}
\end{algorithm}

\section{Computational results}

For comparison, we chose two algorithms presented in this paper and two known heuristic algorithms:
\begin{itemize}
	\item Iterative ILP algorithm from Section~\ref{sec:iterativeILP} (Algorithm~\ref{Alg:ILP});
	
	\item Iterative ILP + LS algorithm from Section~\ref{sec:local_search} (Algorithm~\ref{Alg:ILP_LS});
	
	\item SA: the simulated annealing algorithm that repeatedly finds 2-factors through the reduction to random perfect matching from \cite{kozl:nik:2019};
	
	\item GVNS: the general variable neighborhood search algorithm from \cite{nik:kozl:2020}.
\end{itemize}

The ILP algorithms are implemented in C++, for the SA and GVNS algorithms the existing implementation in Node.js \cite {nik:kozl:2020} is taken. Computational experiments were performed on an Intel (R) Core (TM) i5-4460 machine with a 3.20GHz CPU and 16GB RAM. As the ILP-solver we used SCIP 7.0.2 \cite{gamr:etal:2020ZR}.

The algorithms were tested on random directed and undirected Hamiltonian cycles.
For each graph size, 100 pairs of random permutations with a uniform probability distribution were generated by the Fisher-Yates shuffle algorithm \cite {knuth:1997}.
For two ILP algorithms, a limit of 2 hours was set for each set of 100 instances.
Therefore, the tables indicate how many instances out of 100 the algorithms managed to solve in 2 hours
For both heuristic algorithms, a limit of 60 seconds per test was set, as well as a limit on the number of iterations: $2\,500$ for SA and $250$ for GVNS. The reason is that the heuristic algorithms have a one-sided error. If the algorithm finds a solution, then the solution exists. However, the heuristic algorithms cannot guarantee that the solution to the problem does not exist, only that the solution has not been found in a given time or number of iterations.
For each set of 100 instances, the tables show the average running time in seconds and the average number of iterations separately for feasible and infeasible problems.

{\small
	\begin{table}[p]
		\centering
		\caption{Computational results for pairs of random directed Hamiltonian cycles}
		\label{table:random_directed}
		\begin{tabular}{|*{13}{r|}}
			\hline
			& \multicolumn{6}{c|}{Iterative ILP} & \multicolumn{6}{c|}{Iterative ILP + LS} \\ 
			\hline
			& \multicolumn{3}{c|}{Feasible} & \multicolumn{3}{c|}{Infeasible} & \multicolumn{3}{c|}{Feasible} & \multicolumn{3}{c|}{Infeasible} \\ \hline
			$|V|$ & N & time (s) & Iter & N & time (s) & Iter & N & time (s) & Iter & N & time (s) & Iter \\ 
			\hline
			192 & $21$ & $0.028$ & $4.23$ & $79$ & $0.029$ & $4.22$ & $21$ & $0.022$ & $2.00$ & $79$ & $0.037$ & $3.44$\\
			\hline
			256 & $25$ & $0.111$ & $7.04$ & $75$ & $0.060$ & $5.74$ & $25$ & $0.075$ & $3.12$ & $75$ & $0.090$ & $4.64$\\
			\hline
			384 & $20$ & $0.054$ & $4.65$ & $80$ & $0.082$ & $5.75$ & $20$ & $0.080$ & $2.60$ & $80$ & $0.156$ & $4.33$\\
			\hline
			512 & $22$ & $0.105$ & $5.45$ & $78$ & $0.114$ & $5.58$ & $22$ & $0.125$ & $2.36$ & $78$ & $0.248$ & $4.33$\\
			\hline
			768 & $19$ & $0.148$ & $6.05$ & $81$ & $0.125$ & $5.43$ & $19$ & $0.226$ & $2.26$ & $81$ & $0.439$ & $4.27$\\
			\hline
			1024 & $17$ & $0.182$ & $5.70$ & $83$ & $0.222$ & $6.21$ & $17$ & $0.277$ & $1.88$ & $83$ & $0.925$ & $4.80$\\
			\hline
			1536 & $16$ & $0.325$ & $6.43$ & $84$ & $0.404$ & $6.83$ & $16$ & $1.099$ & $2.50$ & $84$ & $2.312$ & $5.34$\\
			\hline
			2048 & $15$ & $0.568$ & $7.33$ & $85$ & $0.503$ & $6.70$ & $15$ & $3.137$ & $3.13$ & $85$ & $3.829$ & $5.09$\\
			\hline
			3072 & $21$ & $1.130$ & $7.95$ & $79$ & $1.009$ & $7.65$ & $21$ & $4.661$ & $2.42$ & $79$ & $10.722$ & $5.91$\\
			\hline
			4096 & $18$ & $1.681$ & $8.16$ & $82$ & $1.522$ & $7.95$ & $18$ & $18.560$ & $4.44$ & $82$ & $21.283$ & $5.98$\\
			\hline
			& \multicolumn{6}{c|}{SA (perfect matching)} & \multicolumn{6}{c|}{GVNS} \\ 
			\hline
			& \multicolumn{3}{c|}{Solved} & \multicolumn{3}{c|}{Not solved} & \multicolumn{3}{c|}{Solved} & \multicolumn{3}{c|}{Not solved} \\ \hline
			$|V|$ & N & time (s) & Iter & N & time (s) & Iter & N & time (s) & Iter & N & time (s) & Iter \\ 
			\hline
			192 & $21$ & $0.594$ & $183.95$ & $71$ & $8.196$ & $2500$ & $21$ & $0.015$ & $3.76$ & $79$ & $1.618$ & $250$\\
			\hline
			256 & $20$ & $2.557$ & $466.75$ & $80$ & $13.776$ & $2500$ & $25$ & $0.084$ & $10.72$ & $75$ & $2.168$ & $250$\\ 
			\hline
			384 & $17$ & $5.824$ & $500.70$ & $83$ & $28.136$ & $2500$ & $20$ & $0.110$ & $8.00$ & $80$ & $5.958$ & $250$\\
			\hline
			512 & $16$ & $6.525$ & $315.68$ & $84$ & $47.791$ & $2500$ & $22$ & $0.136$ & $5.45$ & $78$ & $10.903$ & $250$\\
			\hline
			768 & $13$ & $13.384$ & $292.15$ & $87$ & $60.000$ & $1420$ & $19$ & $0.643$ & $10.52$ & $81$ & $25.673$ & $250$\\
			\hline
			1024 & $9$ & $11.576$ & $127.88$ & $91$ & $60.000$ & $749.81$ & $17$ & $1.746$ & $16.94$ & $83$ & $42.272$ & $250$\\
			\hline
			1536 & $-$ & $-$ & $-$ & $100$ & $60.000$ & $361.40$ & $16$ & $1.072$ & $5.81$ & $84$ & $60.000$ & $196.64$\\
			\hline
			2048 & $1$ & $8.717$ & $38$ & $99$ & $60.000$ & $249.98$ & $15$ & $3.201$ & $12.00$ & $85$ & $60.000$ & $151.94$\\
			\hline
			3072 & $3$ & $13.719$ & $42.33$ & $97$ & $60.000$ & $185.57$ & $21$ & $5.554$ & $12.28$ & $79$ & $60.000$ & $115.62$\\
			\hline
			4096 & $1$ & $16.713$ & $37$ & $99$ & $60.000$ & $141.15$ & $18$ & $9.395$ & $17.38$ & $82$ & $60.000$ & $98.73$\\
			\hline
		\end{tabular}
	\end{table}
}

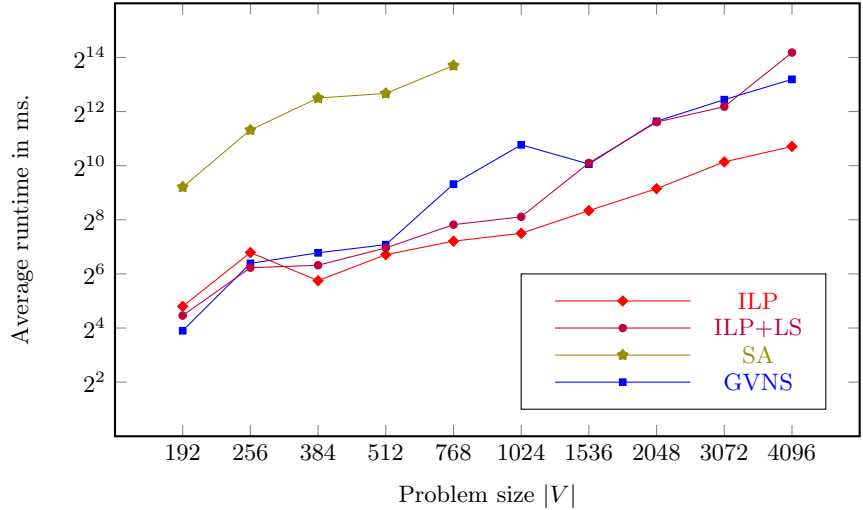
\begin{figure}[p]
	\centering
	\begin{tikzpicture}[scale=1]
	\begin{axis}[
	y=0.36cm,
	x=0.9cm,
	axis line style = thick,
	xlabel={Problem size $|V|$},
	ylabel={Average runtime in ms.},
	xtick       = {1,2,3,4,5,6,7,8,9,10},
	xticklabels = {192,256,384,512,768,1024,1536,2048,3072,4096},
	ytick       = {2,4,6,8,10,12,14},
	yticklabels = {$2^{2}$,$2^{4}$,$2^{6}$,$2^{8}$,$2^{10}$,$2^{12}$,$2^{14}$},
	xmin=0,
	xmax=11,
	ymin=0,
	ymax=16.]

	\node [star,draw,fill,inner sep=1pt,olive] (SA1) at (axis cs: 1,9.21) {};
	\node [star,draw,fill,inner sep=1pt,olive] (SA2) at (axis cs: 2,11.32) {};
	\node [star,draw,fill,inner sep=1pt,olive] (SA3) at (axis cs: 3,12.50) {};
	\node [star,draw,fill,inner sep=1pt,olive] (SA4) at (axis cs: 4,12.67) {};
	\node [star,draw,fill,inner sep=1pt,olive] (SA5) at (axis cs: 5,13.70) {};
	
	\draw [olive] (SA1) -- (SA2) -- (SA3) -- (SA4) -- (SA5);

	\node [draw,fill,inner sep=1.25pt,blue] (GVNS1) at (axis cs: 1,3.90) {};
	\node [draw,fill,inner sep=1.25pt,blue] (GVNS2) at (axis cs: 2,6.39) {};
	\node [draw,fill,inner sep=1.25pt,blue] (GVNS3) at (axis cs: 3,6.78) {};
	\node [draw,fill,inner sep=1.25pt,blue] (GVNS4) at (axis cs: 4,7.08) {};
	\node [draw,fill,inner sep=1.25pt,blue] (GVNS5) at (axis cs: 5,9.32) {};
	\node [draw,fill,inner sep=1.25pt,blue] (GVNS6) at (axis cs: 6,10.77) {};
	\node [draw,fill,inner sep=1.25pt,blue] (GVNS7) at (axis cs: 7,10.06) {};
	\node [draw,fill,inner sep=1.25pt,blue] (GVNS8) at (axis cs: 8,11.64) {};
	\node [draw,fill,inner sep=1.25pt,blue] (GVNS9) at (axis cs: 9,12.44) {};
	\node [draw,fill,inner sep=1.25pt,blue] (GVNS10) at (axis cs: 10,13.19) {};
	
	\draw [blue] (GVNS1) -- (GVNS2) -- (GVNS3) -- (GVNS4) -- (GVNS5) -- (GVNS6) -- (GVNS7) -- (GVNS8) -- (GVNS9) -- (GVNS10);
	
	\node [diamond,draw,fill,inner sep=1pt,red] (ILP1) at (axis cs: 1,4.80) {};
	\node [diamond,draw,fill,inner sep=1pt,red] (ILP2) at (axis cs: 2,6.79) {};
	\node [diamond,draw,fill,inner sep=1pt,red] (ILP3) at (axis cs: 3,5.75) {};
	\node [diamond,draw,fill,inner sep=1pt,red] (ILP4) at (axis cs: 4,6.71) {};
	\node [diamond,draw,fill,inner sep=1pt,red] (ILP5) at (axis cs: 5,7.21) {};
	\node [diamond,draw,fill,inner sep=1pt,red] (ILP6) at (axis cs: 6,7.50) {};
	\node [diamond,draw,fill,inner sep=1pt,red] (ILP7) at (axis cs: 7,8.34) {};
	\node [diamond,draw,fill,inner sep=1pt,red] (ILP8) at (axis cs: 8,9.15) {};
	\node [diamond,draw,fill,inner sep=1pt,red] (ILP9) at (axis cs: 9,10.14) {};
	\node [diamond,draw,fill,inner sep=1pt,red] (ILP10) at (axis cs: 10,10.71) {};
	
	\draw [red] (ILP1) -- (ILP2) -- (ILP3) -- (ILP4) -- (ILP5) -- (ILP6) -- (ILP7) -- (ILP8) -- (ILP9) -- (ILP10);
	
	\node [circle,draw,fill,inner sep=1pt,purple] (ILPLS1) at (axis cs: 1,4.46) {};
	\node [circle,draw,fill,inner sep=1pt,purple] (ILPLS2) at (axis cs: 2,6.23) {};
	\node [circle,draw,fill,inner sep=1pt,purple] (ILPLS3) at (axis cs: 3,6.32) {};
	\node [circle,draw,fill,inner sep=1pt,purple] (ILPLS4) at (axis cs: 4,6.96) {};
	\node [circle,draw,fill,inner sep=1pt,purple] (ILPLS5) at (axis cs: 5,7.82) {};
	\node [circle,draw,fill,inner sep=1pt,purple] (ILPLS6) at (axis cs: 6,8.11) {};
	\node [circle,draw,fill,inner sep=1pt,purple] (ILPLS7) at (axis cs: 7,10.10) {};
	\node [circle,draw,fill,inner sep=1pt,purple] (ILPLS8) at (axis cs: 8,11.61) {};
	\node [circle,draw,fill,inner sep=1pt,purple] (ILPLS9) at (axis cs: 9,12.18) {};
	\node [circle,draw,fill,inner sep=1pt,purple] (ILPLS10) at (axis cs: 10,14.18) {};
	
	\draw [purple] (ILPLS1) -- (ILPLS2) -- (ILPLS3) -- (ILPLS4) -- (ILPLS5) -- (ILPLS6) -- (ILPLS7) -- (ILPLS8) -- (ILPLS9) -- (ILPLS10);

	\draw [red] (axis cs: 6.5,5) -- (axis cs: 8.5,5);
	\node [diamond,draw,fill,inner sep=1pt,red] at (axis cs: 7.5,5) {};
	\node [red] at (axis cs: 9.5,5) {ILP};
	
	\draw [purple] (axis cs: 6.5,4) -- (axis cs: 8.5,4);
	\node [circle,draw,fill,inner sep=1pt,purple] at (axis cs: 7.5,4) {};
	\node [purple] at (axis cs: 9.5,4) {ILP+LS};
	
	\draw [olive] (axis cs: 6.5,3) -- (axis cs: 8.5,3);
	\node [star,draw,fill,inner sep=1pt,olive] at (axis cs: 7.5,3) {};
	\node [olive] at (axis cs: 9.5,3) {SA};
	
	\draw [blue] (axis cs: 6.5,2) -- (axis cs: 8.5,2);
	\node [draw,fill,inner sep=1.25pt,blue] at (axis cs: 7.5,2) {};
	\node [blue] at (axis cs: 9.5,2) {GVNS};
	
	\draw (axis cs: 6,6) -- (axis cs: 10.5,6) -- (axis cs: 10.5,1) -- (axis cs: 6,1) -- cycle;
	
	\end{axis}
	\end{tikzpicture}
	\caption{Computational results for feasible problems on directed graphs}
	\label{image:computational_results_graph_directed}
\end{figure}

The computational results for random directed multigraphs are presented in Table~\ref{table:random_directed} and Figure~\ref{image:computational_results_graph_directed}.
On the considered test set, only 194 instances out of 1\,000 had a solution.
Three algorithms: iterative ILP, iterative ILP + LS, and GVNS correctly solved all instances in the given time, while SA found only 101 Hamiltonian decompositions of 194.
It can be seen that directed multigraphs contain not many subtours.
Thus, the iterative ILP algorithm requires on average only 6.3 iterations to find a solution, and 6.2 iterations to prove that a solution does not exist.
The addition of the local search heuristic to the algorithm makes it possible to reduce the number of iterations by an average of 2.4 times for problems with a solution and 1.3 times for problems without a solution.
In some cases, as for graphs on 192 and 256 vertices, this speeds up the algorithm. However, in most cases, the heuristic does not give an improvement in runtime.
On average, ILP + LS is 3 times slower on problems with a solution and 5 times slower on problems without a solution, and the gap only increases with the growth of the graph size.
We can conclude that a few extra iterations of the ILP-solver turn out to be cheaper in runtime than using an additional heuristic.

Regarding the heuristic algorithms, the performance of GVNS on instances with the existing solution is on average 3.8 times slower than ILP and is comparable to ILP + LS.
While SA completely dropped out of the competition, finding only 101 solutions out of 194, and being an order of magnitude slower.
As for instances without a solution, both heuristic algorithms are not able to determine this scenario and exit only when the limit on the running time or the number of iterations is reached.
In this case, GVNS turns out to be on average 100 times slower than ILP.
However, it is difficult to compare performance here since the time and iteration limits in both heuristic algorithms are set as parameters.
Note that, in the GVNS, the limit was 250 iterations, while the algorithm found a solution, if it exists, on average in 10 iterations. This means that the limit can potentially be lowered to speed up the algorithm. However, this will increase the danger of losing the existing solution.

{\small
\begin{table}[p]
	\centering
	\caption{Computational results for pairs of random undirected Hamiltonian cycles}
	\label{table:random_undirected}
		\begin{tabular}{|*{13}{r|}}
			\hline
			& \multicolumn{6}{c|}{Iterative ILP} & \multicolumn{6}{c|}{Iterative ILP + LS} \\ 
			\hline
			& \multicolumn{3}{c|}{Feasible} & \multicolumn{3}{c|}{Infeasible} & \multicolumn{3}{c|}{Feasible} & \multicolumn{3}{c|}{Infeasible} \\ \hline
			$|V|$ & N & time (s) & Iter & N & time (s) & Iter & N & time (s) & Iter & N & time (s) & Iter \\ 
			\hline
			192 & $100$ & $2.029$ & $23.28$ & $-$ & $-$ & $-$ & $100$ & $0.052$ & $1.24$ & $-$ & $-$ & $-$\\
			\hline
			256 & $100$ & $4.800$ & $30.53$ & $-$ & $-$ & $-$ & $100$ & $0.094$ & $1.30$ & $-$ & $-$ & $-$\\
			\hline
			384 & $100$ & $12.168$ & $34.13$ & $-$ & $-$ & $-$ & $100$ & $0.150$ & $1.27$ & $-$ & $-$ & $-$\\
			\hline
			512 & $100$ & $24.914$ & $44.22$ & $-$ & $-$ & $-$ & $100$ & $0.217$ & $1.28$ & $-$ & $-$ & $-$\\
			\hline
			768 & $100$ & $67.382$ & $54.41$ & $-$ & $-$ & $-$ & $100$ & $0.488$ & $1.29$ & $-$ & $-$ & $-$\\
			\hline
			1024 & $20$ & $396.215$ & $95.40$ & $-$ & $-$ & $-$ & $100$ & $0.721$ & $1.21$ & $-$ & $-$ & $-$\\
			\hline
			1536 & $1$ & $30.598$ & $33$ & $-$ & $-$ & $-$ & $100$ & $1.518$ & $1.34$ & $-$ & $-$ & $-$\\
			\hline
			2048 & $4$ & $1618.87$ & $235.6$ & $-$ & $-$ & $-$ & $100$ & $3.281$ & $1.32$ & $-$ & $-$ & $-$\\
			\hline
			3072 & $4$ & $1772.42$ & $143.25$ & $-$ & $-$ & $-$ & $100$ & $6.746$ & $1.34$ & $-$ & $-$ & $-$\\
			\hline
			4096 & $2$ & $3506.19$ & $168.50$ & $-$ & $-$ & $-$ & $100$ & $14.447$ & $1.38$ & $-$ & $-$ & $-$\\
			\hline
			& \multicolumn{6}{c|}{SA (perfect matching)} & \multicolumn{6}{c|}{GVNS} \\ 
			\hline
			& \multicolumn{3}{c|}{Solved} & \multicolumn{3}{c|}{Not solved} & \multicolumn{3}{c|}{Solved} & \multicolumn{3}{c|}{Not solved} \\ \hline
			$|V|$ & N & time (s) & Iter & N & time (s) & Iter & N & time (s) & Iter & N & time (s) & Iter \\ 
			\hline
			192 & $100$ & $0.884$ & $105.84$ & $-$ & $-$ & $-$ & $100$ & $0.023$ & $1.00$ & $-$ & $-$ & $-$\\
			\hline
			256 & $100$ & $1.904$ & $124.75$ & $-$ & $-$ & $-$ & $100$ & $0.035$ & $1.00$ & $-$ & $-$ & $-$\\ 
			\hline
			384 & $100$ & $7.734$ & $228.22$ & $-$ & $-$ & $-$ & $100$ & $0.073$ & $1.00$ & $-$ & $-$ & $-$\\
			\hline
			512 & $99$ & $12.880$ & $236.39$ & $1$ & $60.000$ & $1016$ & $100$ & $0.133$ & $1.00$ & $-$ & $-$ & $-$\\
			\hline
			768 & $70$ & $21.223$ & $194.74$ & $30$ & $60.000$ & $498.96$ & $100$ & $0.291$ & $1.00$ & $-$ & $-$ & $-$\\
			\hline
			1024 & $46$ & $21.548$ & $124.23$ & $54$ & $60.000$ & $313.14$ & $100$ & $0.511$ & $1.00$ & $-$ & $-$ & $-$\\
			\hline
			1536 & $25$ & $26.140$ & $70.24$ & $75$ & $60.000$ & $157.05$ & $100$ & $1.085$ & $1.00$ & $-$ & $-$ & $-$\\
			\hline
			2048 & $12$ & $35.540$ & $54.33$ & $88$ & $60.000$ & $91.71$ & $100$ & $1.824$ & $1.00$ & $-$ & $-$ & $-$\\
			\hline
			3072 & $6$ & $29.225$ & $19.50$ & $94$ & $60.000$ & $41.11$ & $100$ & $4.235$ & $1.00$ & $-$ & $-$ & $-$\\
			\hline
			4096 & $-$ & $-$ & $-$ & $100$ & $60.000$ & $22.22$ & $100$ & $7.593$ & $1.00$ & $-$ & $-$ & $-$\\
			\hline
		\end{tabular}
\end{table}
}

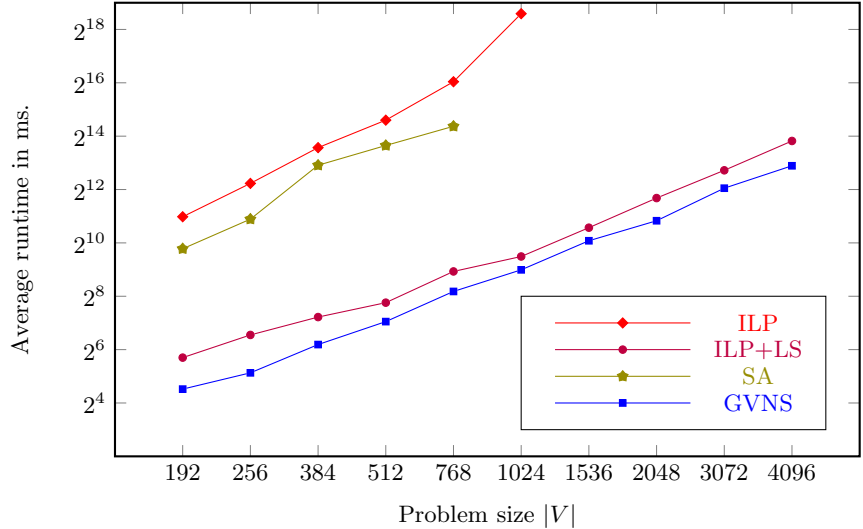
\begin{figure}[p]
	\centering
	\begin{tikzpicture}[scale=1]
	\begin{axis}[
	y=0.355cm,
	x=0.9cm,
	axis line style = thick,
	xlabel={Problem size $|V|$},
	ylabel={Average runtime in ms.},
	xtick       = {1,2,3,4,5,6,7,8,9,10},
	xticklabels = {192,256,384,512,768,1024,1536,2048,3072,4096},
	ytick       = {4,6,8,10,12,14,16,18},
	yticklabels = {$2^{4}$,$2^{6}$,$2^{8}$,$2^{10}$,$2^{12}$,$2^{14}$,$2^{16}$,$2^{18}$},
	xmin=0,
	xmax=11,
	ymin=2,
	ymax=19.]

	\node [star,draw,fill,inner sep=1pt,olive] (SA1) at (axis cs: 1,9.78) {};
	\node [star,draw,fill,inner sep=1pt,olive] (SA2) at (axis cs: 2,10.89) {};
	\node [star,draw,fill,inner sep=1pt,olive] (SA3) at (axis cs: 3,12.91) {};
	\node [star,draw,fill,inner sep=1pt,olive] (SA4) at (axis cs: 4,13.65) {};
	\node [star,draw,fill,inner sep=1pt,olive] (SA5) at (axis cs: 5,14.37) {};
	
	\draw [olive] (SA1) -- (SA2) -- (SA3) -- (SA4) -- (SA5);

	\node [draw,fill,inner sep=1.25pt,blue] (GVNS1) at (axis cs: 1,4.52) {};
	\node [draw,fill,inner sep=1.25pt,blue] (GVNS2) at (axis cs: 2,5.13) {};
	\node [draw,fill,inner sep=1.25pt,blue] (GVNS3) at (axis cs: 3,6.19) {};
	\node [draw,fill,inner sep=1.25pt,blue] (GVNS4) at (axis cs: 4,7.05) {};
	\node [draw,fill,inner sep=1.25pt,blue] (GVNS5) at (axis cs: 5,8.18) {};
	\node [draw,fill,inner sep=1.25pt,blue] (GVNS6) at (axis cs: 6,8.99) {};
	\node [draw,fill,inner sep=1.25pt,blue] (GVNS7) at (axis cs: 7,10.08) {};
	\node [draw,fill,inner sep=1.25pt,blue] (GVNS8) at (axis cs: 8,10.83) {};
	\node [draw,fill,inner sep=1.25pt,blue] (GVNS9) at (axis cs: 9,12.05) {};
	\node [draw,fill,inner sep=1.25pt,blue] (GVNS10) at (axis cs: 10,12.89) {};
	
	\draw [blue] (GVNS1) -- (GVNS2) -- (GVNS3) -- (GVNS4) -- (GVNS5) -- (GVNS6) -- (GVNS7) -- (GVNS8) -- (GVNS9) -- (GVNS10);

	\node [diamond,draw,fill,inner sep=1pt,red] (ILP1) at (axis cs: 1,10.98) {};
	\node [diamond,draw,fill,inner sep=1pt,red] (ILP2) at (axis cs: 2,12.23) {};
	\node [diamond,draw,fill,inner sep=1pt,red] (ILP3) at (axis cs: 3,13.57) {};
	\node [diamond,draw,fill,inner sep=1pt,red] (ILP4) at (axis cs: 4,14.60) {};
	\node [diamond,draw,fill,inner sep=1pt,red] (ILP5) at (axis cs: 5,16.04) {};
	\node [diamond,draw,fill,inner sep=1pt,red] (ILP6) at (axis cs: 6,18.59) {};
	
	\draw [red] (ILP1) -- (ILP2) -- (ILP3) -- (ILP4) -- (ILP5) -- (ILP6);
	
	\node [circle,draw,fill,inner sep=1pt,purple] (ILPLS1) at (axis cs: 1,5.7) {};
	\node [circle,draw,fill,inner sep=1pt,purple] (ILPLS2) at (axis cs: 2,6.55) {};
	\node [circle,draw,fill,inner sep=1pt,purple] (ILPLS3) at (axis cs: 3,7.22) {};
	\node [circle,draw,fill,inner sep=1pt,purple] (ILPLS4) at (axis cs: 4,7.76) {};
	\node [circle,draw,fill,inner sep=1pt,purple] (ILPLS5) at (axis cs: 5,8.93) {};
	\node [circle,draw,fill,inner sep=1pt,purple] (ILPLS6) at (axis cs: 6,9.49) {};
	\node [circle,draw,fill,inner sep=1pt,purple] (ILPLS7) at (axis cs: 7,10.57) {};
	\node [circle,draw,fill,inner sep=1pt,purple] (ILPLS8) at (axis cs: 8,11.68) {};
	\node [circle,draw,fill,inner sep=1pt,purple] (ILPLS9) at (axis cs: 9,12.72) {};
	\node [circle,draw,fill,inner sep=1pt,purple] (ILPLS10) at (axis cs: 10,13.82) {};
	
	\draw [purple] (ILPLS1) -- (ILPLS2) -- (ILPLS3) -- (ILPLS4) -- (ILPLS5) -- (ILPLS6) -- (ILPLS7) -- (ILPLS8) -- (ILPLS9) -- (ILPLS10);

	\draw [red] (axis cs: 6.5,7) -- (axis cs: 8.5,7);
	\node [diamond,draw,fill,inner sep=1pt,red] at (axis cs: 7.5,7) {};
	\node [red] at (axis cs: 9.5,7) {ILP};
	
	\draw [purple] (axis cs: 6.5,6) -- (axis cs: 8.5,6);
	\node [circle,draw,fill,inner sep=1pt,purple] at (axis cs: 7.5,6) {};
	\node [purple] at (axis cs: 9.5,6) {ILP+LS};
	
	\draw [olive] (axis cs: 6.5,5) -- (axis cs: 8.5,5);
	\node [star,draw,fill,inner sep=1pt,olive] at (axis cs: 7.5,5) {};
	\node [olive] at (axis cs: 9.5,5) {SA};
	
	\draw [blue] (axis cs: 6.5,4) -- (axis cs: 8.5,4);
	\node [draw,fill,inner sep=1.25pt,blue] at (axis cs: 7.5,4) {};
	\node [blue] at (axis cs: 9.5,4) {GVNS};
	
	\draw (axis cs: 6,8) -- (axis cs: 10.5,8) -- (axis cs: 10.5,3) -- (axis cs: 6,3) -- cycle;
	
	\end{axis}
	\end{tikzpicture}
	\caption{Computational results for undirected graphs}
	\label{image:computational_results_graph_undirected}
\end{figure}

The situation for undirected multigraphs is fundamentally different. 
It is known that random undirected regular graphs have a Hamiltonian decomposition with a very high probability, which allows finding the decomposition asymptotically almost surely by random matchings in polynomial time \cite {kim:2001}. 
This approach is in some way similar to the considered SA algorithm.
In our case, the problem is slightly different, since the multigraph $x \cup y$ always has a decomposition into cycles $x$ and $y$, and we need to find another decomposition into cycles $z$ and $w$.
Nevertheless, for all $1\,000$ instances on undirected cycles (Table~\ref{table:random_undirected} and Figure~\ref{image:computational_results_graph_undirected}), there was a second Hamiltonian decomposition, and the vertices of the $\operatorname{TSP}(n)$ polytope were not adjacent. From a geometric point of view, this means that the degrees of vertices in 1-skeleton are much less than the total number of vertices, so two random vertices are not adjacent with a very high probability.

Summary for random undirected multigraphs: both iterative ILP + LS and GVNS solved all 1\,000 instances, SA solved 558 instances, and ILP solved only 531 instances in a given time.
It can be concluded that the iterative ILP algorithm was not very successful for undirected graphs and showed similar results to the SA algorithm. On instances up to 768 vertices, where all tests were solved, the ILP was on average 2.3 times slower than the SA.
The problem is that undirected multigraphs contain a large number of subtours that have to be forbidden. On average, the ILP algorithm took about 86 iterations to find a solution.
On the other hand, the addition of the local search heuristic to the ILP algorithm reduced the running time by an average of 200 times, and the number of iterations by 65 times.
The ILP + LS algorithm showed results similar to GVNS, solving all test instances and being on average only 1.8 times slower.
This time loss is due to two factors.
Firstly, the GVNS has a more complex heuristic with several neighborhood structures, which made it possible to find all solutions in just 1 iteration.
Secondly, one iteration of the ILP-solver is much more expensive than constructing 2-factors through the reduction to perfect matching.

It should be noted that although all 1\,000 random instances on undirected graphs had a solution, in the general case, the traveling salesperson polytope contains adjacent vertices for which, accordingly, the Hamiltonian decomposition does not exist. Moreover, the 1-skeleton of the traveling salesperson polytope has cliques with an exponential number of vertices \cite{bond:1983}.
Thus, the ILP + LS algorithm may turn out to be more promising, since it will be able to prove that there is no Hamiltonian decomposition for the given problem.

{
	\begin{table}[t]
		\centering
		\caption{Computational results for infeasible problems on undirected pyramidal tours}
		\label{table:random_undirected_pyramidal}
		\begin{tabular}{|*{9}{r|}}
			\hline
			& \multicolumn{2}{c|}{Iterative ILP} & \multicolumn{2}{c|}{Iterative ILP+LS} & \multicolumn{2}{c|}{SA} & \multicolumn{2}{c|}{GVNS} \\ 
			\hline
			$|V|$ & time (s) & Iter & time (s) & Iter & time (s) & Iter & time (s) & Iter  \\ 
			\hline
			128 & $0.177$ & $10.90$ & $0.239$ & $8.94$ & $3.439$ & $2500$ & $1.495$ & $250$ \\
			\hline
			192 & $0.418$ & $14.78$ & $0.645$ & $11.26$ & $6.901$ & $2500$ & $2.278$ & $250$ \\
			\hline
			256 & $0.822$ & $19.40$ & $1.575$ & $14.84$ & $11.687$ & $2500$ & $3.253$ & $250$ \\
			\hline
			384 & $2.425$ & $28.50$ & $5.549$ & $20.59$ & $25.215$ & $2500$ & $5.609$ & $250$ \\
			\hline
			512 & $4.586$ & $36.08$ & $13.488$ & $26.24$ & $43.904$ & $2500$ & $9.868$ & $250$ \\
			\hline
			768 & $17.740$ & $56.41$ & $56.824$ & $41.25$ & $60.000$ & $1713$ & $17.313$ & $250$ \\
			\hline
		\end{tabular}
	\end{table}
}

We ran additional tests to investigate this scenario (Table~\ref{table:random_undirected_pyramidal}). 
Using the vertex adjacency criterion for the pyramidal tours polytope \cite{bond:nik:2018}, we generated 6 groups of 50 pairs of such undirected pyramidal tours $x$ and $y$ that the multigraph $x \cup y$ is guaranteed not to contain a Hamiltonian decomposition into cycles $z$ and $w$.
It can be seen that although the additional local search heuristic reduced the number of iterations by an average of 1.3 times, the total running time increased by an average of 2.2 times.
Indeed, the local search takes extra time to find a solution that does not exist.
However, this slight slowdown is acceptable, given that on undirected multigraphs with an existing solution (Table~\ref{table:random_undirected}), the local search heuristic gives an average speed up of 200 times.
Note that for undirected graphs the number of iterations grows significantly faster than for directed multigraphs (Tables~\ref{table:random_directed}) since the undirected multigraphs contain a large number of subtours that have to be forbidden.
Nevertheless, the ILP algorithms have the advantage over the SA and GVNS here, since the heuristic algorithms cannot guarantee that the problem is infeasible.

\section{Conclusion}

We introduced two iterative ILP algorithms to find a Hamiltonian decomposition of the 4-regular multigraph.
On random undirected multigraphs, the version enhanced by the local search heuristic turned out to be much more efficient than the basic ILP algorithm, showing results comparable to the known general variable neighborhood search heuristic.
While for random directed multigraphs the iterative ILP algorithm significantly surpassed in speed the previously known algorithms.
The key feature that distinguishes the ILP algorithms from previously known heuristics is that they can prove that the Hamiltonian decomposition in the graph does not exist.

The directions for further development are as follows.
Firstly, we can consider a more complex heuristic with several neighborhood structures, as in \cite{nik:kozl:2020}, to speed up the algorithm on problems with an existing solution.
Secondly, it is of great interest to add to the model other classes of facet inequalities of the traveling salesperson polytope, like
2-matching and clique-tree inequalities \cite{grot:padb:1985}, that can significantly reduce the number of expensive calls of an ILP-solver.

\subsubsection*{Acknowledgements.}
	We are very grateful to the anonymous reviewers for their comments and suggestions which helped to improve the
	presentation of the results in this paper.

%
%
%
%

\end{document}